\documentclass{article}
\usepackage{amsmath, amssymb, graphicx, tikz}
\usetikzlibrary{shapes.geometric}
\usetikzlibrary{arrows, automata}
\usepackage{caption}
\usepackage{cite}
\usepackage{url}
\usepackage{hyperref}
\usepackage{verbatim}

\newtheorem{theorem}{Theorem}[section]
\newtheorem{corollary}{Corollary}[theorem]
\newtheorem{lemma}[theorem]{Lemma}
\newtheorem{example}{Example}[section]

\title{Nash Equilibria of Rock Paper Scissors Variants}
\author{Adrian Thananopavarn}
\date{September 2024}

\begin{document}

\maketitle
\begin{abstract}
We generalize Rock Paper Scissors to complete directed graphs, or tournaments, on $n$ vertices. Properties of the mixed-strategy Nash equilibria of these tournaments are discussed, particularly those with Nash equilibria where all of the strategies have a nonzero probability. We find graph-theoretic properties of such games and tabulate them for $n \leq 7$.
\end{abstract}

\tableofcontents
\pagebreak
\section{Introduction}
Rock Paper Scissors is a classic 2-player game where each player simultaneously chooses either Rock, Paper, or Scissors to display using a hand gesture, and the winner is determined using a directed graph $G$:

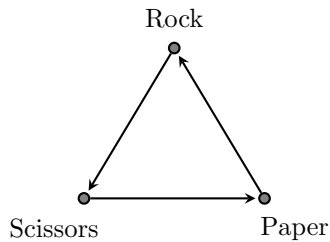
\begin{figure}[!htbp]
\hspace*{4cm}
\begin{tikzpicture}[thick,scale=0.8, > = stealth, 
            shorten > = 1pt, 
            auto]
\tikzstyle{every node}=[circle, draw, fill=black!50, inner sep=0pt, minimum width=4pt]
\node[draw=white, fill=white] at (5.5,1.5) {Paper};
\node[draw=white, fill=white] at (3.5,5) {Rock};
\node[draw=white, fill=white] at (1.5,1.5) {Scissors};
\node[] (v0) at (2,2) {};
\node[] (v1) at (5,2) {};
\node[] (v2) at (3.5,4.5) {};
\draw[->] (v1)--(v2);
\draw[->] (v0)--(v1);
\draw[->] (v2)--(v0);

\end{tikzpicture}
\caption{A directed graph $G$ demonstrating which player wins in Rock Paper Scissors. If there exists an edge pointing from Player 1's choice to Player 2's choice, then Player 1 wins.}
\end{figure}

Let a \textit{strategy} in such a game be one of the options to play in the game. Rock Paper Scissors, as it is usually played, has only three strategies: Rock, Paper, and Scissors. In this paper, we discuss generalizations to the game, with the number of strategies $n > 3$.

We only concern ourselves with games where for any two distinct strategies $s_1, s_2$, either $s_1$ beats $s_2$ or $s_2$ beats $s_1$. The result is a complete directed graph, or \textit{tournament}.

Such games have been studied before for special cases. Akin \cite{akin} investigates Eulerian tournaments, where each vertex has the same number of edges entering as leaving. Using the automorphism groups of the tournament, Akin classifies Eulerian tournaments and finds all Eulerian tournaments up to isomorphism with $n \leq 9$.

For $n \leq 5$, such games have been studied in the more general case, where the graph need not be directed \cite{inproceedings}, by Tagiew. If two vertices do not have an edge between them, their respective strategies tie against each other. Tagiew finds 4 games out of the 582 possible directed graphs that have unique mixed strategy Nash equilibria where all strategies have positive probability. 2 of these games are complete directed graphs, which are also found in this paper.

In general, though, the properties of Nash equilibria of such tournaments are as yet unexplored. In section 2, we investigate the Nash equilibria of Rock Paper Scissors variants, finding that they are unique (Theorem 2.6) and have an odd number of nonzero entries (Theorem 2.8). In section 3, we present several constructions of variants with known properties. In section 4, we use these properties, connecting Nash equilibria to the graph-theoretic properties of the tournament. In section 5, we tabulate several Rock-Paper Scissors variants with $n \leq 7$.

\section{Nash equilibria}

In this section, we discuss the game-theoretic properties of Rock Paper Scissors variants, which tell us which strategies one should play for a given variant. In games such as these, often no single strategy is the best, because each strategy is beaten by at least one other strategy. So, it is helpful to define a notion for picking each strategy with a certain probability. A \textit{strategy profile} $\vec{a} = \langle a_1, a_2, \dots, a_n \rangle$ is an $n$-dimensional vector where each entry $a_k$ is the probability that one plays strategy $k$. By probability axioms, the sum of all $a_k$ is $1$. We may evaluate the relative quality of two independent strategy profiles $\vec{a}, \vec{b}$ by comparing the probabilities that $\vec{a}$ beats $\vec{b}$,  $\vec{a}$ draws with $\vec{b}$, and $\vec{a}$ loses to $\vec{b}$. These will be respectively denoted as $P(W)$, $P(D)$, and $P(L)$.

A mixed-strategy Nash equilibrium (henceforth referred to as just a \textit{Nash equilibrium}) is a strategy profile $\vec{a}$ such that for every strategy profile $\vec{b}$, $P(W) \geq P(L)$. A pivotal result in game theory \cite{nash} is stated below without proof.

\begin{theorem}[Nash Equilibria, Nash 1950]
All games with a finite number of players and finitely many strategies have Nash equilibria.
\end{theorem}

A quick application of this shows that all Rock Paper Scissors variants (with $n$ finite) have a Nash equilibrium. The key issue is finding a specific Nash equilibrium for an arbitrary Rock Paper Scissors variant $G$.

\subsection{Uniqueness of Nash Equilibria}

Here, we prove the uniqueness of Nash equilibria for Rock Paper Scissors variants by representing the conditions for Nash equilibria with a system of inequalities. We have $P(W) + P(D) + P(L) = 1$, which gives an equivalence:

$$P(W) \geq P(L) \Longleftrightarrow 2P(W) + P(D) \geq 1.$$
\vspace*{0cm}

Then, the quantity $2P(W) + P(D)$ may be represented as follows:
\begin{align*}
2P(W) + P(D) &= \sum_{i=1}^{n}\sum_{j=1}^{n}g_{ij}a_ib_j\\
&= \sum_{j=1}^{n}b_j\sum_{i=1}^{n}g_{ij}a_i
\end{align*}

where $g_{ij}$ is the outcome of the game: $2$ if strategy $i$ beats $j$, $1$ if $i=j$, and $0$ if strategy $i$ loses to $j$. Notably, $g_{ij} + g_{ji} = 2$, because every pair of distinct strategies has a winner and a loser. This motivates the following lemma, which gives a better characterization of Nash equilibria of Rock Paper Scissors variants.

\begin{lemma} A strategy profile $\vec{a}$ is a Nash equilibrium iff $\sum_{i=1}^{n}g_{ij}a_i \geq 1$ for $1 \leq j \leq n$.
\end{lemma}

\textit{Proof}: Let 
$$c_j(\vec{a}) = \sum_{i=1}^{n}g_{ij}a_i.$$
for $1 \leq j \leq n$. Then, from the above equation, we have 
$$2P(W) + P(D) = \sum_{j=1}^{n}c_j(\vec{a})b_j$$.

Now, if $\vec{a}$ is a Nash equilibrium, then for every strategy profile $\vec{b}$, 
$$\sum_{j=1}^{n}c_j(\vec{a})b_j \geq 1.$$
If any particular $c_k(\vec{a})$ were less than $1$, one could let the corresponding $b_k = 1$, and set all other $b_j$ equal to $0$. Then, we would have $\sum_{j=1}^{n}c_j(\vec{a})b_j = c_k(\vec{a}) < 1.$ So, all $c_k(\vec{a})$ must be at least $1$.

Conversely, if all $c_j(\vec{a}) \geq 1$, then for any $\vec{b}$, 

$$\sum_{j=1}^{n}c_j(\vec{a})b_j \geq \sum_{j=1}^{n}b_j = 1.$$

Therefore, $\vec{a}$ is a Nash equilibrium. $\square$

\vspace*{0.5cm}
Intuitively, $c_j(\vec{a})$ is the expected outcome of the game given that one uses strategy profile $\vec{a}$ and the opponent always chooses strategy $j$. For a Nash equilibrium, the above lemma gives a lower bound of $c_j(\vec{a}) \geq 1$. The next lemma shows what happens if $c_j(\vec{a}) > 1$.

\begin{lemma} If $\vec{a}$ and $\vec{b}$ are both Nash equilibria, and $c_j(\vec{a}) > 1$ for some $j$ with $1 \leq j \leq n$, then $b_j = 0$.
\end{lemma}

\textit{Proof}: Suppose this were not the case: without loss of generality, assume that $c_n(\vec{a}) > 1$ but $b_n > 0$. Because $\vec{a}$ and $\vec{b}$ are both Nash equilibria, we have: 

\begin{align*}
&\hspace{0.7cm} P(W) = P(L) \\
&\Rightarrow 2P(W) + P(D) = 1 \\
&\Rightarrow \sum_{j=1}^{n}c_j(\vec{a})b_j = 1 \\
&\Rightarrow 1 = c_n(\vec{a})b_n + \sum_{j=1}^{n-1}c_j(\vec{a})b_j > b_n + \sum_{j=1}^{n-1}b_j = 1.\\
\end{align*}
By contradiction, $b_n$ must equal $0$. $\square$

\vspace*{0.5cm}
Letting $\vec{a} = \vec{b}$, we immediately get the following corollary:

\begin{corollary}
If $\vec{a}$ is a Nash equilibrium and $c_j(\vec{a}) > 1$ for some $j$ with $1 \leq j \leq n$, then $a_j = 0$.
\end{corollary}

By contraposition of this corollary, if all $a_j > 0$ in a Nash equilibrium, $c_j(\vec{a}) = 1$, so we have a system of equations:

$$\sum_{i=1}^{n}g_{ij}a_i = 1 \hspace{1cm} (1 \leq j \leq n)$$

Expressed as a matrix:

\[
    \begin{bmatrix} 
    g_{11} & \dots  & g_{n1}\\
    \vdots & \ddots & \vdots\\
    g_{1n} & \dots  & g_{1n} 
    \end{bmatrix}
    \vec{a}
    = 
    \begin{bmatrix} 
    1 \\
    \vdots \\
    1
    \end{bmatrix}
\]

In this matrix, from our construction of the $g_{ij}$, the main diagonal is all $1$, and the rest of the entries are all even. Let this matrix be denoted $G$, and the vector of all 1's denoted as $\vec{1}$. This leads to the next lemmas, which find that Nash equilibria in these games are unique.

\begin{lemma}
If $\vec{a}$ is a Nash equilibrium and $a_j > 0$ for all $1 \leq j \leq n$, then $\vec{a}$ is the only Nash equilibrium.
\end{lemma}

\textit{Proof}: First, notice that if $a_j > 0$ for $1 \leq j \leq n$, then all Nash equilibria $\vec{b}$ have $c_j(\vec{b}) = 1$, by Lemma $2.3$.

Now, all that remains to show is that $G\vec{b} = \vec{1}$ has a unique solution.

To compute the determinant of $G$ \cite{dets}:

$$\text{det}(G) = \sum_{\sigma \in S_n} \left( \text{sgn}(\sigma) \prod_{k=1}^n g_{k\sigma(k)}\right)$$

where $\sigma$ is an element of the permutation group $S_n$, and $\text{sgn}(\sigma)$ is the signature of $\sigma$. We will disregard the signature (which is always $\pm 1$) and focus on the even-odd parity of the expression. All non-identity permutations $\sigma$ will have some element $g_{k\sigma(k)}$ such that $k \neq \sigma(k)$, which is even, but the identity permutation has all $g_{kk} = 1$. So, 

$$\text{det}(G) = \pm \prod_{k=1}^n g_{kk} + \sum_{\sigma \in S_n \setminus I} \left( \pm \prod_{k=1}^n g_{k\sigma(k)}\right)$$

is an odd number plus the sum of even numbers, which is nonzero. Therefore, $G\vec{b} = \vec{1}$ has a unique solution. Since $\vec{a}$ is a solution, it is the only solution, so it is the only Nash equilibrium. $\square$

\begin{lemma}
If $\vec{a}$ is a Nash equilibrium and $c_j(\vec{a}) = 1$ for all $1 \leq j \leq n$, then every $a_j$, when represented as a fraction in lowest terms $p_j/q_j$, has both $p_j$ and $q_j$ odd.
\end{lemma}

\textit{Proof}: 
We first show that all $q_j$ are odd.

Suppose this were not the case: that some $a_k$ had an even denominator when written in lowest terms. Then, express each fraction with a common denominator: let $q = \text{lcm}(q_1,\dots,q_n)$, and let $p'_j = p_jq/q_j$. Now, as $a_j = p'_j/q$,

$$\sum_{i=1}^{n}g_{ij}p'_j = q \hspace{1cm} (1 \leq j \leq n)$$

As $q$ is even by our assumption, and $g_{ij}$ is even whenever $i \neq j$, $g_{jj}p'_j = p'_j$ must be even for all $j$. However, the $q_k$ with the highest power of $2$ will have $q/q_k$ odd, so $p_k$ would have to be even to make $p'_k$ even. This contradicts that $p_k/q_k$ is in lowest terms. So, all $q_j$ are odd.

\vspace*{0.5cm}

Now, to show all $p_j$ are odd, we use the same equation $\sum_{i=1}^{n}g_{ij}p'_j = q$, except now we know that $q$ is odd. So, since $g_{ij}$ is even whenever $i \neq j$, $g_{jj}p'_j = p'_j$ must be odd for all $j$. So, all $p_j$ are odd. $\square$

\begin{corollary}
If $\vec{a}$ is a Nash equilibrium and $c_j(\vec{a}) = 1$ for all $1 \leq j \leq n$, then $a_j > 0$ for all $1 \leq j \leq n$.
\end{corollary}

Now, using Lemma $2.4$, we get this corollary:

\begin{corollary}
If $\vec{a}$ is a Nash equilibrium and $c_j(\vec{a}) = 1$ for all $1 \leq j \leq n$, then $\vec{a}$ is the only Nash equilibrium.
\end{corollary}

Finally, the following theorem shows that Nash equilibria are unique for general $G$.

\begin{theorem}[Uniqueness of Nash Equilibria]
    For any variant $G$, there is exactly one Nash equilibrium $\vec{a}$.
\end{theorem}

\textit{Proof}: Let $\vec{a}$ be a Nash equilibrium. We wish to show that it is unique.

Suppose that $m$ of the $c_j(\vec{a})$ are equal to $1$. Without loss of generality, let $c_j(\vec{a}) = 1$ for $1 \leq j \leq m$, and $c_j(\vec{a}) > 1$ for $m < j \leq n$.

By Corollary $2.3.1$, $a_j = 0$ for $m < j \leq n$. Now, if we let $G_m$ be the matrix of the first $m$ rows and $m$ columns of $G$, and let $\vec{a}_m$ be the first $m$ entries of $\vec{a}$, we have $G_m\vec{a}_m = [1]$. $G_m$ is the outcome matrix of the game only consisting of the first $m$ strategies of the original game, and using Corollary $2.5.1$, $a_j > 0$ for $1 < j \leq m$.

Now, consider any Nash equilibrium $\vec{b}$. Again by Corollary $2.3.1$, $b_j = 0$ for $m < j \leq n$, and $c_j(\vec{b}) = 1$ for $1 \leq j \leq m$. We then have $G_m\vec{b_m} = [1]$. However, by Corollary $2.5.2$, $\vec{a_m}$ was the only solution to this equation, so $\vec{a}_m = \vec{b}_m$.

Combining these with the zero entries from $m+1$ to $n$, we get $\vec{a}=\vec{b}$, so $\vec{a}$ is the unique Nash equilibrium. $\square$

\vspace*{0.5cm}

This result allows us to henceforth refer to the Nash equilibrium of a Rock Paper Scissors game uniquely.

\subsection{Parity of Nash Equilibria}

For any tournament on $n$ vertices, we have exactly one Nash equilibrium. However, the Nash equilibrium could have some, or many, strategies with probability $0$. From a game design perspective, this results in a less varied strategy profile and a less interesting game.

The following theorem quantifies this:
\begin{theorem}[Reduction of Strategies with 0 Probability]
Let $\vec{a}$ be the Nash equilibrium for a game $G$ where $a_j > 0$ for $1 \leq j \leq m$ and $a_j = 0$ for $1 < m \leq n$. Suppose also that $\vec{b}$ is the Nash equilibrium for the same game $G_m$ with only the first $m$ strategies. Then, $\vec{b}$ = $\vec{a}_m$.
\end{theorem}

\textit{Proof}: For $1 \leq j \leq m$, we have

$$\sum_{i=1}^{n}g_{ij}a_i = c_j(\vec{a}) = 1$$.

However, as $a_j = 0$ for $1 < m \leq n$, we may omit all but the first $m$ terms of the sum:

$$\sum_{i=1}^{m}g_{ij}a_i = 1$$.

Now, since for all $i,j \leq m$, the matrix $G$ is the same as $G_m$, $\vec{a}_m$ is a Nash equilibrium for $G_m$. So, $\vec{b}$ = $\vec{a}_m$ by uniqueness. $\square$

\vspace*{0.5cm}

Let a Rock-Paper-Scissors variant be \textit{all-positive} if every probability in its Nash equilibrium is positive. The above theorem tells us that the strategies with zero probability in the Nash equilibrium do not affect the Nash equilibrium for the strategies with positive probability. Any game that has strategies with zero probability in the Nash equilibrium, then, may be reduced to an all-positive game. This paper will primarily investigate all-positive games from this point onwards.

\begin{theorem}[Odd Parity of All-Positive Games]
    All all-positive games have odd $n$.
\end{theorem}
\textit{Proof:} Suppose we have an all-positive game with matrix $G$ and Nash equilibrium $\vec{a}$. Similar to the proof of Lemma 2.5, we may represent all $a_j$
as $p_j/q_j$ in lowest terms, or $p'_j/q$ using the least common denominator. From Lemma 2.4, all $p'_j$ and $q$ are odd.

As probabilities sum to 1,

$$\sum_{j=1}^n p'_j = q.$$

Thus, there must be an odd number of terms in the sum, so $n$ must be odd. $\square$

This result may explain why most published Rock Paper Scissors variants have an odd number of strategies \cite{umop}!

\subsection{Computation of Nash Equilibria in All-Positive games}

Given an all-positive game, the Nash equilibrium is the solution to a system of equations. But given any game, how may we determine whether the game is all-positive (and thus may be solved using a system of equations)? As an example, we compute the Nash equilibrium in Rock Paper Scissors.

\begin{example}[Nash equilibrium of Rock Paper Scissors]
\end{example}

Assuming that Rock Paper Scissors is all-positive, the equation for Rock Paper Scissors is as follows:

\[
    \begin{bmatrix} 
    1 & 2 & 0\\
    0 & 1 & 2\\
    2 & 0 & 1 
    \end{bmatrix}
    \vec{a}
    = 
    \begin{bmatrix} 
    1 \\
    1 \\
    1
    \end{bmatrix}
\]

Using row-reduction, the unique solution is $\vec{a} = \langle 1/3, 1/3, 1/3 \rangle$. Now, checking this strategy profile, it satisfies both the probability axioms and the Nash equilibrium inequalities, so this strategy profile is the Nash equilibrium.

\vspace*{0.5cm}

In general, if one solves the system of linear equalities, then the solution will automatically satisfy the Nash equilibrium inequalities. The strategy profile will then be a Nash equilibrium if it satisfies the probability axioms: namely that all probabilities are non-negative and that they sum to $1$.

The following lemma shows that one of these conditions is automatic: the solution will sum to $1$.

\begin{lemma}
For any variant with $n$ odd, the solution to $G \vec{a} = [1]$ will have a sum of $1$.
\end{lemma}

\textit{Proof}: We consider the following matrix equation:

\[
    \begin{bmatrix} 
    g_{11} - 1 & \dots  & g_{n1} - 1\\
    \vdots & \ddots & \vdots\\
    g_{1n} - 1 & \dots  & g_{1n} - 1\\
    \end{bmatrix}
    \vec{a}
    = 
    \begin{bmatrix} 
    0 \\
    \vdots \\
    0 \\ 
    \end{bmatrix}
\]

Here, all elements of the original matrix $G$ are decreased by 1, to make a new matrix $G'$. This resulting matrix is anti-symmetric, meaning $-G' = (G')^T$.

To find the determinant of $G'$, we again use the formula:

$$\text{det}(G) = \sum_{\sigma \in S_n} \left( \text{sgn}(\sigma) \prod_{k=1}^n g_{k\sigma(k)}'\right).$$

Now, every permutation $\sigma \in S_n$ has a unique inverse, so:

$$2 * \text{det}(G) = \sum_{\sigma \in S_n} \left( \text{sgn}(\sigma) \prod_{k=1}^n g_{k\sigma(k)}' +  \text{sgn}(\sigma^{-1}) \prod_{k=1}^n g_{k\sigma^{-1}(k)}'\right)$$

$$= \sum_{\sigma \in S_n} \left( \text{sgn}(\sigma) \prod_{k=1}^n g_{k\sigma(k)}' +  \text{sgn}(\sigma) \prod_{k=1}^n g_{\sigma(k)k}'\right)$$

$$= \sum_{\sigma \in S_n} \left( \text{sgn}(\sigma) \prod_{k=1}^n g_{k\sigma(k)}' + (-1)^n\text{sgn}(\sigma) \prod_{k=1}^n g_{k\sigma(k)}'\right) = 0.$$

Therefore, at least one of the rows of $G'$ is a linear combination of the others. Replace this row with a row of only $1$'s to make $G''$:

\[
    \begin{bmatrix} 
    g_{11} - 1 & \dots  & g_{n1} - 1\\
    \vdots & \vdots & \vdots\\
    1 & 1 & 1\\
    \vdots & \vdots & \vdots\\
    g_{1n} - 1 & \dots  & g_{1n} - 1\\
    \end{bmatrix}
    \vec{a}
    = 
    \begin{bmatrix} 
    0 \\
    \vdots \\
    1 \\
    \vdots \\
    0 \\ 
    \end{bmatrix}
\]

Now, to solve this new matrix equation, we again find the determinant. To ease this process, we add the row of $1$'s to each other row, which does not change the determinant. In the interest of not having a triple prime, let this new matrix be $H$.  

$$\text{det}(H) = \sum_{\sigma \in S_n} \left( \text{sgn}(\sigma) \prod_{k=1}^n h_{k\sigma(k)}\right).$$

Similarly to the proof of Lemma 2.4, the only permutation that only contains odd values of $h_{k\sigma(k)}$ is the identity permutation, so det($H$) is odd. Therefore, $G''$ is invertible, and our matrix equation has a unique solution $\vec{a}$.

Notably, $\vec{a}$ is also a solution to our original equation $G\vec{a} = [1]$, so it is the only solution. Furthermore, its entries sum to $1$ by the additional row we added to the matrix. Thus, the solution to the original equation has a sum of $1$. $\square$

\vspace*{0.5cm}

Because of this lemma, in order to check if a game with an odd number of strategies is all-positive, one may solve the matrix equation; if all of the entries are positive, the game is all-positive and the solution is the Nash equilibrium, and if not all of the entries are positive, the game is not all-positive (and the solution need not be the Nash equilibrium).

\section{Families of Rock Paper Scissors Variants}
There are several infinite families of Rock Paper Scissors variants that we have studied, two of which are all-positive, and one of which is not.

\subsection{Eulerian Tournaments}

A tournament is \textit{Eulerian} if every strategy beats exactly half of the other strategies \cite{akin}. Because there must be an even number of other strategies, all Eulerian tournaments have odd $n$. Up to isomorphism, the normal Rock Paper Scissors game is the only Eulerian tournament of size $n = 3$. 

As an example, an Eulerian Rock Paper Scissors variant with $n = 5$ is shown below. This variant is called Rock Paper Scissors Lizard Spock \cite{kass}.

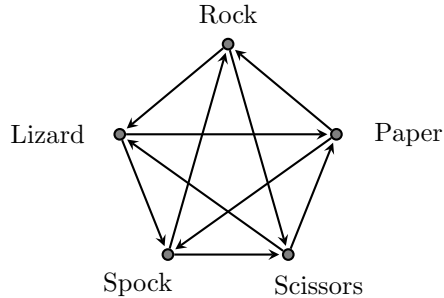
\begin{figure}[!htbp]
\hspace*{3cm}
\begin{tikzpicture}[thick,scale=0.8, > = stealth, shorten > = 1pt, auto]
\tikzstyle{every node}=[circle, draw, fill=black!50, inner sep=0pt, minimum width=4pt]
\node[draw=white, fill=white] at (8,6) {Paper};
\node[draw=white, fill=white] at (5,8) {Rock};
\node[draw=white, fill=white] at (2,6) {Lizard};
\node[draw=white, fill=white] at (3.5,3.5) {Spock};
\node[draw=white, fill=white] at (6.5,3.5) {Scissors};
\node[] (v0) at (3.2,6) {};
\node[] (v1) at (4,4) {};
\node[] (v2) at (6,4) {};
\node[] (v3) at (6.8,6) {};
\node[] (v4) at (5,7.5) {};
\draw[->] (v4)--(v0);
\draw[->] (v0)--(v1);
\draw[->] (v1)--(v2);
\draw[->] (v2)--(v3);
\draw[->] (v3)--(v4);
\draw[->] (v4)--(v2);
\draw[->] (v2)--(v0);
\draw[->] (v0)--(v3);
\draw[->] (v3)--(v1);
\draw[->] (v1)--(v4);

\end{tikzpicture}
\caption{Rock Paper Scissors Lizard Spock, an Eulerian variant with $n=5$.}
\end{figure}

\begin{theorem}
In any Eulerian game, the Nash equilibrium is to choose each strategy with probability $1/n$.
\end{theorem}

\textit{Proof}: To show that the strategy profile with all entries $a_i = 1/n$ satisfies the equations, we must ensure that:

$$\sum_{i=1}^{n}g_{ij}a_i = 1 \hspace{1cm} (1 \leq j \leq n)$$

In any Eulerian game, we have:

$$\sum_{i=1}^{n}g_{ij}/n = \frac{1}{n} \sum_{i=1}^{n}g_{ij} = \frac{1}{n} \left( 1*1 + 0* \frac{n-1}{2} + 2* \frac{n-1}{2}\right) = 1.$$

So, $\langle 1/n, 1/n, \dots, 1/n \rangle$ satisfies the equations, so it is the Nash equilibrium. $\square$

\vspace*{0.5cm}
Eulerian tournaments exist for every odd $n$, and examples have been constructed for $n$ as large as $101$ \cite{umop}. 

\subsection{Game Substitution}
For any two variants $A$ and $B$, one may replace any strategy $a$ in $A$ with the entire game $B$ to make a new game $C$ with all strategies in $A \setminus a$ and $B$, as follows:

\begin{itemize}
    \item If strategy $a_1$ beat $a_2$ in $A$, and neither of $a_1$ or $a_2$ are $a$, then $a_1$ still beats $a_2$ in $C$.
    \item If strategy $b_1$ beat $b_2$ in $A$, then $b_1$ still beats $b_2$ in $C$.
    \item Iff strategy $a_1$ beat $a$ in $A$, then $a_1$ beats all strategies of $B$ in $C$.
\end{itemize}

For example, Figure 3 shows an example of a game of $5$ strategies, where one of the strategies in ordinary Rock Paper Scissors is substituted with another Rock Paper Scissors game.

\begin{figure}[!htbp]
\hspace*{4cm}
\begin{tikzpicture}[thick,scale=0.8, > = stealth, shorten > = 1pt, auto]
\tikzstyle{every node}=[circle, draw, fill=black!50, inner sep=0pt, minimum width=4pt]
\node[draw=white, fill=white] at (2.75,7) {$1$};
\node[draw=white, fill=white] at (4.5,7.75) {$2$};
\node[draw=white, fill=white] at (6,7) {$3$};
\node[draw=white, fill=white] at (7,4) {$5$};
\node[draw=white, fill=white] at (2,4) {$4$};
\node[] (v0) at (4.5,7.25) {};
\node[] (v1) at (2.5,4) {};
\node[] (v2) at (6.5,4) {};
\node[] (v3) at (5.5,6.75) {};
\node[] (v4) at (3.5,6.75) {};
\draw[->] (v0)--(v4);
\draw[->] (v3)--(v0);
\draw[->] (v4)--(v3);
\draw[->] (v1)--(v4);
\draw[->] (v1)--(v0);
\draw[->] (v1)--(v3);
\draw[->] (v4)--(v2);
\draw[->] (v0)--(v2);
\draw[->] (v3)--(v2);
\draw[->] (v2)--(v1);

\end{tikzpicture}
\caption{A Rock Paper Scissors variant with $n=5$. This game may be seen as Rock Paper Scissors with $n=3$, where the strategies are to pick strategy 4, 5, or any of $\{1,2,3\}$. 1, 2, and 3 form another Rock Paper Scissors game, which substitutes for a strategy in the first. We call this game Rock Paper Scissors Turtle Balloon.}
\end{figure}
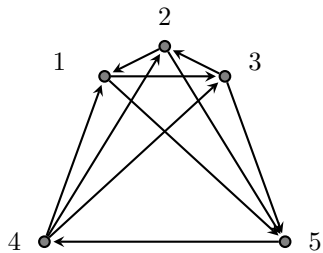

If one knows the Nash equilibria for games $A$ and $B$, this determines the Nash equilibrium for $A$ with the game $B$ substituted for strategy $a \in A$.

\begin{theorem}
The Nash equilibrium for a game $A$ with $m$ strategies, where strategy $m$ has been substituted with a game $B$ with $n$ strategies, may be determined using the Nash equilibria for the component games, $\vec{a}$ and $\vec{b}$. The Nash equilibrium for the new game is $\langle a_1, a_2, \dots ,a_{m-1}, a_mb_1, a_mb_2,\dots,a_mb_n \rangle$.
\end{theorem}

\textit{Proof}: Let the matrices for games $A$ and $B$ be $G$ and $H$, respectively. As we have Nash equilibria for games $A$ and $B$, we know that:

$$\sum_{i=1}^{m}g_{ij}a_i \geq 1,  1 \leq j \leq m-1$$

$$\sum_{i=1}^{n}h_{ij}b_i \geq 1,  1 \leq j \leq n-1.$$

Now, to check that these inequalities hold for the new game with our proposed Nash equilibrium, we construct a new game matrix for the new game, according to the substitution rules:

\[
\left[
\begin{array}{c c c|c c c}
g_{11} & \dots & g_{m-1,1} & g_{m1} & \dots & g_{m1}\\
\vdots & \ddots & \vdots & \vdots & \cdots & \vdots\\
g_{1,m-1} & \dots & g_{m-1,m-1} & g_{m,m-1} & \dots & g_{m,m-1}\\
\hline
g_{1m} & \dots & g_{m-1,m} & h_{11} & \dots & h_{n1}\\
\vdots & \vdots & \vdots & \vdots & \ddots & \vdots\\
g_{1m} & \dots & g_{m-1,m} & h_{1n} & \dots & h_{n1}\\
\end{array}
\right]
\]
\vspace*{0.5cm}

For the upper section of this matrix,
$$\sum_{i=1}^{m-1}g_{ij}a_i + \sum_{i=1}^{n}g_{mj}a_mb_i = \sum_{i=1}^{m-1}g_{ij}a_i + a_m\sum_{i=1}^{n}g_{mj}b_i \geq \sum_{i=1}^{m-1}g_{ij}a_i + g_{mj}a_m \geq 1.$$
For the lower section,
$$\sum_{i=1}^{m-1}g_{im}a_i + \sum_{i=1}^{n}h_{ij}a_mb_i = \sum_{i=1}^{m-1}g_{im}a_i + a_m\sum_{i=1}^{n}h_{ij}b_i \geq \sum_{i=1}^{m-1}g_{im}a_i + a_m \geq 1$$.

As all of the inequalities hold, $\langle a_1, a_2, \dots ,a_{m-1}, a_mb_1, a_mb_2,\dots,a_mb_n \rangle$ is the Nash equilibrium for the new variant. $\square$

\subsection{Dominating Strategies}
Suppose we have the following game:

\begin{figure}[!htbp]
\hspace*{4cm}
\begin{tikzpicture}[thick,scale=0.8, > = stealth, shorten > = 1pt, auto]
\tikzstyle{every node}=[circle, draw, fill=black!50, inner sep=0pt, minimum width=4pt]
\node[draw=white, fill=white] at (2.5,7.5) {Rock};
\node[draw=white, fill=white] at (6.5,7.5) {Scissors 1};
\node[draw=white, fill=white] at (6.5,3.5) {Scissors 2};
\node[draw=white, fill=white] at (2.5,3.5) {Paper};
\node[] (v0) at (3,7) {};
\node[] (v1) at (3,4) {};
\node[] (v2) at (6,4) {};
\node[] (v3) at (6,7) {};
\draw[->] (v0)--(v3);
\draw[->] (v3)--(v1);
\draw[->] (v1)--(v0);
\draw[->] (v0)--(v2);
\draw[->] (v2)--(v1);
\draw[->] (v3)--(v2);

\end{tikzpicture}
\caption{An example game that functions the same as Rock Paper Scissors, but with two Scissors options. Scissors 1 beats Scissors 2.}
\end{figure}
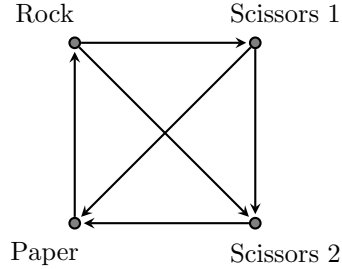

The Nash equilibrium for this game, which may be verified by checking the Nash equilibrium inequalities, is to choose each of Rock, Paper, and Scissors 1 $1/3$ of the time, and to never choose Scissors 2. Thus, the game is not all-positive. This example motivates the following definition: let strategy $i$ \textit{dominate} strategy $j$ if $i$ beats $j$, and there exist no other strategies $k$ such that $j$ beats $k$ and $k$ beats $i$. In this example, Scissors 1 dominates Scissors 2.

\begin{theorem}
    If strategy $1$ dominates strategy $2$, then the Nash equilibrium $\vec{a}$ has $a_2 = 0$.
\end{theorem}

\textit{Proof}: Suppose $a_2$ were nonzero. Then, consider the distinct strategy profile $\langle a_1+a_2, 0, a_3, a_4,\dots,a_n \rangle$. We shall show that this is also a Nash equilibrium. To do this, we show that

$$g_{1j}(a_1+a_2) + \sum_{i=3}^{n}g_{ij}a_i \geq \sum_{i=1}^{n}g_{ij}a_i \geq 1,$$

where the second inequality follows from the fact that $\vec{a}$ is a Nash equilibrium.

The difference in the first inequality is 
$$g_{1j}(a_1+a_2) + \sum_{i=3}^{n}g_{ij}a_i - \sum_{i=1}^{n}g_{ij}a_i = a_2(g_{1j} - g_{2j}).$$

For $j=1$, this quantity is positive because $g_{11} = 1, g_{21} = 0$.

For $j=2$, this quantity is positive because $g_{12} = 2, g_{22} = 1$.

For $3 \leq j \leq n$, both $g_{1j}$ and $g_{2j}$ are either $0$ or $2$. However, the case where $g_{1j}=0$, $g_{2j}=2$ cannot occur because this would mean that 2 beats $j$ and $j$ beats 1. So, the difference is non-negative.

Therefore, $\langle a_1+a_2, 0, a_3, a_4,\dots,a_n \rangle$ is a Nash equilibrium, contradicting uniqueness. So, $a_2$ must equal $0$. $\square$

\section{Graph-Theoretic Properties}

Theorem 3.3 gives a graph-theoretic way to automatically show that a tournament does not yield an all-positive game. In this section, we shall investigate properties of games where no strategy dominates another.

\subsection{King Chickens}
A 1980 paper by Maurer titled \textit{The King Chicken Theorems} \cite{chicken} is about tournaments on chickens; a chicken $a$ \textit{pecks} another chicken $b$ if there is an edge from $a$ to $b$. A chicken $a$ is a \textit{king chicken} if for all other chickens $b$, either $a$ pecks $b$, or there exists a chicken $c$ such that $a$ pecks $c$ and $c$ pecks $b$.

If we let strategies in our game be chickens, then strategy $a$ being a king chicken is equivalent to $a$ not being dominated by any other strategy. This means that if no strategy dominates another, then all strategies are king chickens. Let a game where no strategy dominates another (i.e. all strategies are king chickens) be called a \textit{royal flock}. All all-positive games are royal flocks, but not all royal flocks are necessarily all-positive. This may be seen from the following theorem.

\begin{theorem} [Royal Flocks, Maurer 1980]
There exist royal flocks of all positive integer sizes except $n=2$ and $n=4$.
\end{theorem}

By this theorem, there exists a royal flock of size $6$, but as this game has an even number of strategies, it is not all-positive. This game is shown below.

\begin{figure}[!htbp]
\hspace*{4cm}
\begin{tikzpicture}[thick,scale=0.8, > = stealth, shorten > = 1pt, auto]
\tikzstyle{every node}=[circle, draw, fill=black!50, inner sep=0pt, minimum width=4pt]
\node[draw=white, fill=white] at (5,3.5) {$c_4$};
\node[draw=white, fill=white] at (7.5,4.75) {$c_3$};
\node[draw=white, fill=white] at (7.5,7.25) {$c_2$};
\node[draw=white, fill=white] at (5,8.5) {$c_1$};
\node[draw=white, fill=white] at (2.5,7.25) {$c_6$};
\node[draw=white, fill=white] at (2.5,4.75) {$c_5$};
\node[] (v0) at (5,8) {};
\node[] (v1) at (3,7) {};
\node[] (v2) at (3,5) {};
\node[] (v3) at (5,4) {};
\node[] (v4) at (7,5) {};
\node[] (v5) at (7,7) {};
\draw[->] (v0)--(v1);
\draw[->] (v1)--(v2);
\draw[->] (v2)--(v3);
\draw[->] (v3)--(v4);
\draw[->] (v4)--(v5);
\draw[->] (v5)--(v0);
\draw[->] (v1)--(v5);
\draw[->] (v5)--(v3);
\draw[->] (v3)--(v1);
\draw[->] (v2)--(v4);
\draw[->] (v4)--(v0);
\draw[->] (v0)--(v2);
\draw[->] (v0)--(v3);
\draw[->] (v4)--(v1);
\draw[->] (v2)--(v5);

\end{tikzpicture}
\caption{Maurer's example of a royal flock with $n=6$. It is not all-positive.}
\end{figure}
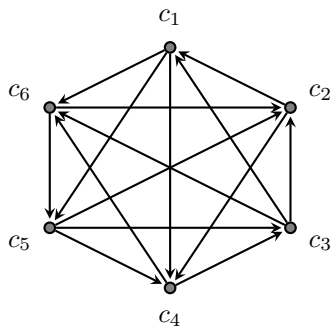

\subsection{Known Edges in Royal Flocks}
When searching for all-positive games, the following theorems are helpful for providing a few starting edges on the graph.

\begin{theorem}
All royal flocks of size $n\geq3$ have an $n$-cycle.
\end{theorem}
Proof: Suppose for the sake of contradiction that $G$ is a royal flock of size $n \geq 3$ without an $n$-cycle. Then, let $k$ be the size of the largest cycle in $G$. We know a cycle must exist: if we choose any two distinct vertices $a$ and $b$ where $a$ beats $b$, there must be another vertex $c$ where $b$ beats $c$ and $c$ beats $a$, by our definition of king chicken. This makes a cycle with $a$, $b$, and $c$.

Let $v_1, v_2, \dots, v_k, v_{k+1}=v_1$ be the vertices of the largest cycle, where $v_1$ beats $v_2$, $v_2$ beats $v_3$, and so on until $v_k$ beats $v_1$. Then, as $k<n$, let $w$ be a vertex not in this cycle.

Suppose $w$ beats any vertex in the cycle; without loss of generality let this vertex be $v_1$. Then, if $w$ does not beat all other $v$ in the cycle, let $v_j$ be the last $v$ in the cycle such that $v_j$ beats $w$. Then, we have a $k+1$-cycle, $v_1, v_2, \dots, v_j, w, v_{j+1}, \dots, v_k, v_{k+1}=v_1$,  which is a contradiction.

If $w$ does beat all $v$ in the cycle, $w$ cannot dominate $v_1$. So, there must be another vertex $u$ such that $v_1$ beats $u$, and $u$ beats $w$. Then, we have a cycle $v_1, u, w, v_2,\dots, v_k, v_{k+1}=v_1$ of size $k+2$, again a contradiction.

If $w$ does not beat any vertex in the cycle, $v_2$ cannot dominate $w$. So, there must be another vertex $u$ such that $w$ beats $u$ and $u$ beats $v_2$. Then, we have a cycle $v_1, w, u, v_2,\dots, v_k, v_{k+1}=v_1$ of size $k+2$, again a contradiction.

In any case, we have constructed a larger cycle in $G$, which contradicts the assumption that the largest cycle was size $k$. So, all royal flocks of size $n\geq3$ have an $n$-cycle. $\square$

\vspace*{0.5cm}

\begin{theorem}
   No royal flock has a strategy with in-degree or out-degree $n-1$. If a royal flock has a strategy with in-degree or out-degree $n-2$, the game is Rock Paper Scissors, but with one strategy substituted for a game of size $n-2$.
\end{theorem}

\textit{Proof}: If a tournament has a strategy with in-degree $n-1$, then it is dominated by all other strategies. If a tournament has a strategy with out-degree $n-1$, then it dominates all other strategies. In either case, the tournament is not a royal flock.

If a royal flock has a strategy $a$ with in-degree $n-2$, then for each of the $n-2$ strategies $b$ that beat $a$, there must be a strategy $c$ where $a$ beats $c$ and $c$ beats $b$. $c$ is unique, as $a$ only beats $1$ other strategy. So, $a$ beats $c$, $c$ beats all other strategies, and all other strategies beat $a$. Similarly, if a strategy $c$ has out-degree $n-2$, then there must be a strategy $a$ where $a$ beats $c$, $c$ beats all other strategies, and all other strategies beat $a$. In either case, the tournament is Rock Paper Scissors, but with one strategy substituted for a game of size $n-2$. $\square$

\section{Prime Rock Paper Scissors Variants}

Let a \textit{prime} Rock Paper Scissors variant be a variant that cannot be represented as the substitution of one game inside another. In this section, we compile a list of all-positive prime Rock Paper Scissors variants with size $n \leq 7$, found using a mixture of graph-theoretic and exhaustive techniques. We name the variants $n_k$, where $n$ is the number of strategies and $k$ is a ranking of the game by magnitude of the Nash equilibrium vector.

In the interest of drawing fewer edges, in this section, to tell if one vertex beats another, use the following rule: if there is no edge between them, the higher vertex wins, and if there is an edge between them, the lower vertex wins. The vertices are indexed top to bottom. Moreover, to simplify the Nash equilibria, all Nash equilibrium vectors will be multiplied by their least common denominator, to make all values integers.

\subsection{n = 1, 3 ,5}
If $n=1$, then there is trivially only one variant. This variant is like Rock Paper Scissors, if one was only allowed to pick Rock.

\vspace*{0.5cm}
If $n=3$, then there must be a $3$-cycle, by Theorem 4.2. This uniquely determines the game as regular Rock Paper Scissors. 

\vspace*{0.5cm}
If $n=5$, there must be a $5$-cycle, by Theorem 4.2. By Theorem 4.3, for the game to be prime, all vertices must have in-degree 2 and out-degree 2, making the game Eulerian. So, the 5 edges not in the first cycle must form a cycle of their own. There is only one graph satisfying these conditions up to isomorphism: Rock Paper Scissors Lizard Spock.

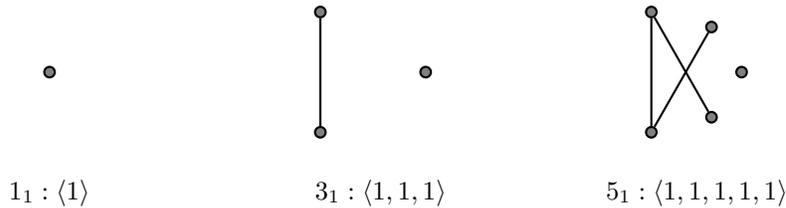
\begin{figure}[!htbp]
\hspace*{1cm}
\begin{tikzpicture}[thick,scale=0.8, auto]
\tikzstyle{every node}=[circle, draw, fill=black!50, inner sep=0pt, minimum width=4pt]
\node[draw=white, fill=white] at (3.5,4) {$1_1: \langle 1 \rangle$};
\node[draw=white, fill=white] at (9,4) {$3_1: \langle 1,1,1 \rangle$};
\node[draw=white, fill=white] at (14.25,4) {$5_1: \langle 1,1,1,1,1 \rangle$};
\node[] (v0) at (3.5,6) {};
\node[] (v1) at (8,7) {};
\node[] (v2) at (8,5) {};
\node[] (v3) at (9.75,6) {};
\node[] (v4) at (13.5,7) {};
\node[] (v5) at (13.5,5) {};
\node[] (v6) at (14.5,6.75) {};
\node[] (v7) at (15,6) {};
\node[] (v8) at (14.5,5.25) {};
\draw (v1)--(v2);
\draw (v4)--(v5);
\draw (v4)--(v8);
\draw (v5)--(v6);

\end{tikzpicture}
\end{figure}

\subsection{n = 7}
The following graphs were generated using Mathematica code from Appendix A. This code generates all possible adjacency matrices for tournaments that satisfy Theorem 4.2. The generated graphs were then tested for isomorphism, and $8$ prime variants were found.

\begin{figure}[!htbp]

\begin{tikzpicture}[thick,scale=0.61, > = stealth , auto]
\tikzstyle{every node}=[circle, draw, fill=black!50, inner sep=0pt, minimum width=4pt]
\node[draw=white, fill=white] at (2.5,13.7) {$7_1: \langle 1,1,1,1,1,1,1 \rangle$};
\node[draw=white, fill=white] at (7.25,12.5) {$7_2: \langle 1,1,1,1,1,1,1 \rangle$};
\node[draw=white, fill=white] at (12,13.7) {$7_3: \langle 1,1,1,1,1,1,1 \rangle$};
\node[draw=white, fill=white] at (16.75,12.5) {$7_4: \langle 7,5,1,9,1,5,7 \rangle$};
\node[draw=white, fill=white] at (2.5,6.7) {$7_5: \langle 3,1,1,1,1,1,3 \rangle$};
\node[draw=white, fill=white] at (7.25,5.5) {$7_6: \langle 5,3,1,7,3,1,5 \rangle$};
\node[draw=white, fill=white] at (12,6.7) {$7_7: \langle 3,5,1,1,1,5,3 \rangle$};
\node[draw=white, fill=white] at (16.75,5.5) {$7_8: \langle 7,3,1,1,1,3,7 \rangle$};
\node[] (v0) at (1.5,18.5) {};
\node[] (v1) at (1.5,14.5) {};
\node[] (v2) at (2.75,14.75) {};
\node[] (v3) at (4,15.5) {};
\node[] (v4) at (4.5,16.5) {};
\node[] (v5) at (4,17.5) {};
\node[] (v6) at (2.75,18.25) {};
\node[] (v7) at (6.25,18.5) {};
\node[] (v8) at (7.5,18.25) {};
\node[] (v9) at (8.75,17.5) {};
\node[] (v10) at (9.25,16.5) {};
\node[] (v11) at (8.75,15.5) {};
\node[] (v12) at (7.5,14.75) {};
\node[] (v13) at (6.25,14.5) {};
\node[] (v14) at (11,18.5) {};
\node[] (v15) at (12.25,18.25) {};
\node[] (v16) at (13.5,17.5) {};
\node[] (v17) at (14,16.5) {};
\node[] (v18) at (13.5,15.5) {};
\node[] (v19) at (12.25,14.75) {};
\node[] (v20) at (11,14.5) {};
\node[] (v21) at (15.75,18.5) {};
\node[] (v22) at (17,18.25) {};
\node[] (v23) at (18.25,17.5) {};
\node[] (v24) at (18.75,16.5) {};
\node[] (v25) at (17,14.75) {};
\node[] (v26) at (15.75,14.5) {};
\node[] (v27) at (1.5,11.5) {};
\node[] (v28) at (1.5,7.5) {};
\node[] (v29) at (4.5,9.5) {};
\node[] (v30) at (4,10.5) {};
\node[] (v31) at (4,8.5) {};
\node[] (v32) at (18.25,15.5) {};
\node[] (v33) at (2.75,11.25) {};
\node[] (v34) at (2.75,7.75) {};
\node[] (v35) at (6.25,11.5) {};
\node[] (v36) at (6.25,7.5) {};
\node[] (v37) at (7.5,7.75) {};
\node[] (v38) at (7.5,11.25) {};
\node[] (v39) at (8.75,10.5) {};
\node[] (v40) at (8.75,8.5) {};
\node[] (v41) at (11,11.5) {};
\node[] (v42) at (9.25,9.5) {};
\node[] (v43) at (11,7.5) {};
\node[] (v44) at (12.25,7.75) {};
\node[] (v45) at (13.5,8.5) {};
\node[] (v46) at (14,9.5) {};
\node[] (v47) at (13.5,10.5) {};
\node[] (v48) at (12.25,11.25) {};
\node[] (v49) at (15.75,11.5) {};
\node[] (v50) at (17,11.25) {};
\node[] (v51) at (18.25,10.5) {};
\node[] (v52) at (18.75,9.5) {};
\node[] (v53) at (18.25,8.5) {};
\node[] (v54) at (17,7.75) {};
\node[] (v55) at (15.75,7.5) {};
\draw (v0)--(v1);
\draw (v1)--(v6);
\draw (v6)--(v2);
\draw (v2)--(v0);
\draw (v0)--(v3);
\draw (v5)--(v1);
\draw (v7)--(v13);
\draw (v13)--(v8);
\draw (v8)--(v12);
\draw (v12)--(v7);
\draw (v7)--(v10);
\draw (v10)--(v13);
\draw (v9)--(v11);
\draw (v14)--(v20);
\draw (v20)--(v15);
\draw (v15)--(v18);
\draw (v14)--(v19);
\draw (v19)--(v16);
\draw (v14)--(v17);
\draw (v17)--(v20);
\draw (v21)--(v25);
\draw (v25)--(v23);
\draw (v26)--(v22);
\draw (v22)--(v32);
\draw (v21)--(v24);
\draw (v24)--(v26);
\draw (v27)--(v28);
\draw (v28)--(v29);
\draw (v29)--(v27);
\draw (v33)--(v34);
\draw (v34)--(v30);
\draw (v33)--(v31);
\draw (v35)--(v40);
\draw (v38)--(v42);
\draw (v39)--(v36);
\draw (v42)--(v37);
\draw (v35)--(v42);
\draw (v42)--(v36);
\draw (v41)--(v43);
\draw (v41)--(v44);
\draw (v44)--(v46);
\draw (v43)--(v48);
\draw (v48)--(v46);
\draw (v47)--(v45);
\draw (v49)--(v55);
\draw (v55)--(v52);
\draw (v52)--(v49);
\draw (v50)--(v54);
\draw (v51)--(v53);

\end{tikzpicture}
\end{figure}

\section{Future Work}
One could find more Rock-Paper-Scissors variants with 9 or more strategies. There are Eulerian variants of any odd size \cite{akin}, and these may be substituted within each other to create infinitely many variants. We have not discovered any other infinite families of all-positive variants, which could be the subject of future work.

Another potential topic of research could be to find more graph-theoretic properties of these variants: are there stronger conditions on the graphs than royal flock-ness that are required for the variant to be all-positive?

Finally, although we have a general method for finding the Nash equilibrium of all-positive games, what if the game isn't all-positive? Is there a simple way to find the Nash equilibrium using the matrix or graph structure?

\section{Acknowledgements}
Thank you to the students and staff of MathILy 2024, where I presented on this topic and received positive feedback, inspiring the writing of this paper. I especially thank Mathes Miller-Priddy for finding four of the variants with $7$ strategies using an exhaustive search: namely, $7_4$ through $7_7$.

\appendix

\section{Appendix: Code to Generate Variants with n=7}
Written by Mathes Miller-Priddy.

\begin{verbatim}
dim = 7;

allGames = {};
For[n = 0, n < 2^((dim - 2)*(dim - 1)/2), n++,
  Game = Map[Append[#, 1] &, IdentityMatrix[dim]];
  fill = 0;
  For[i = 1, i < dim, i ++,
   For[j = i + 2, j < dim + 1, j++,
     fill++;
     Game[[i, j]] = 2*Quotient[Mod[n, 2^fill], 2^(fill - 1)];
     ];
   ];
  
  For[i = 2, i < dim + 1, i++,
   For[j = 1, j < i + 1, j++,
     Game[[i, j]] = 2 - Game[[j, i]];
     ];
   ];
  
  sol = RowReduce[Game][[All, dim + 1]];
  isPos = True;
  For[i = 1, i <= dim, i++,
   If[sol[[i]] <= 0,
     isPos = False;
     ];	
   ];
  
  If[isPos, 
   allGames = Append[allGames, {Game, sol}];
   ];
  ];
Print[allGames // MatrixForm]
\end{verbatim}

\bibliographystyle{plainurl}
\bibliography{refs}

\end{document}